
\input amstex
\documentstyle{amsppt}
\document

\input amstex
\documentstyle{amsppt}
\magnification=1200
\hsize=16truecm
\tolerance=2000
\vsize=21truecm
\baselineskip 14pt
\topmatter
\title On The Surjectivity Of Wahl Maps \\  On A
General Curve
\endtitle
\author Roberto Paoletti
\endauthor
\address Mathematics Department
UCLA, Los Angeles CA 90024 \endaddress
\endtopmatter
\document
\heading Introduction \endheading
Consider a smooth projective curve $C$ and two line bundles $L$ and $N$ on it.
It is well known that there is a linear map, given by section multiplication
$$\mu: H^0(C,L)\otimes H^0(C,N)\longrightarrow H^0(C,L+N).$$
We define the {\it module of relations} of $L$ and $N$, denoted $R(L,N)$,
to be the kernel of $\mu$. The Wahl map, or Gaussian map
$$\gamma_{L,N}: R(L,N)\longrightarrow H^0(K+L+N)$$
(where $K$ denotes the canonical line bundle on $C$) is defined by making
sense of the expression $\gamma_{L,N}(s,t)=:sdt-tds$.
These maps have attracted increasing attention since Wahl's basic observation
that they relate to the deformation theory of the projective cone over $C$
(\cite {W1}).
In fact, if $L$ is a very ample line bundle on $C$ the cokernels of
$\gamma_{K,L^{i-1}}$, for $i$ positive, are dual to the first order
deformations of the projective cone which smooth the vertex.
{}From this it follows, for example, that if $C$ is the hyperplane section of
a (projective) K3 surface, then $\gamma_{K,K}$ is not surjective.
This was proved from a deformation theoretic point of view by Wahl, and along
different lines by Beauville and Merindol(\cite{BM}).

This circle of ideas has led to the question of the behavior of $\gamma_K
=\gamma_{K,K}$ on a general curve. In fact, this map being onto implies that
the general canonical curve is not the hyperplane section of a smooth
surface. Ciliberto
{\it et al} have proved, using degeneration methods, that for a curve with
general moduli $\gamma_K$ surjects
(\cite{CHM}). Mukai then observed (\cite {M}) that
if $C$ is a smooth curve lying on a K3 surface and such that
the class of $C$ genererates $Pic(C)$, then
on $C$ there are minimal pencils for which the adjoint line bundle
is not projectively normal.
Voisin has then generalized this observation into a new
conceptual approach to the probem
(\cite {V}). Namely, she shows that if $C$ is a Petri
general curve for which $\gamma_K$ is not onto, then on $C$ there exist
complete linear series of dimension one and minimal degree such that for all of
these the
adjoint line bundle is not linearly normal.
In other words, if $A$ is a minimal pencil and $\gamma_K$
is not surjective the multiplication map
$S^2H^0(C,K-A)\longrightarrow H^0(C,2K-2A)$
can't be onto. Then she shows that this cannot
happen on a general curve, thereby proving surjectivity of $\gamma_K$ in this
case.
Her proof uses two very different arguments in the odd genus case and in the
even genus one.

At the same time, there has been growing activity concerning the problem
of the surjectivity of $\gamma _{L,N}$, for arbitrary line bundles
$L$ and $N$ on $C$.
This more general question has been explored, among others,
by Bertram, Ein and Lazarsfeld (\cite{BEL}), and by Wahl (\cite{W1},
\cite{W2}, \cite{W3}, \cite{W4}).
The first three authors have found conditions on the degree
of $L$ involving the Clifford index of the curve that guarantee
surjectivity of $\gamma _{K,L}$. Wahl has found other conditions, and he
formulated a conjecture to the effect
that $\gamma _{K,L}$ is onto as soon as $deg(L)\ge 2g+$ some suitable
constant. More generally, he has posed the question of finding a
geometric interpretation of the failure of $\gamma _{L,N}$ to surject,
and of the resulting stratification of the Picard group of $C$ in terms
of the corank of $\gamma _{L,N}$.
The object of this article is to show that even for these more general
Wahl maps one can still interpret the failure of $\gamma _{L,N}$ to surject
in terms of the existence of pencils of small degree for which
suitable section multiplication maps are not surjective.

Before describing the results, let me recall
that $W^r_d(C)$ denotes the subvariety of $Pic^d(C)$ consisting of the
line bundles $L$ on $C$ satisfying $deg(L)=d$ and $h^0(L)\ge r+1$. Then
the main theorem is

\proclaim {Theorem A} Let $C$ be a general curve of genus $g
>8$ and $L$ be a line
bundle on it. Then
\roster
\item If
$g=2s$ and $deg(L)>
3s$, or if $g=2s+1$ and $deg(L)>
3s+3$ and $\gamma_{K,L}$ is not onto, then section multiplication
$H^0(K-A)\otimes H^0(L-A)\longrightarrow H^0(K+L-2A)$ is not surjective,
for general $A\in W^1_{s+2}$.
\item Suppose that $L$ is chosen generally and that $deg(L)\ge 3s+10$
when $g=2s$ or that $deg(L)\ge 3s+6$ when $g=2s+1$. Then
the above multiplication is onto, for a general choice of such an $A$.
\endroster
\endproclaim


{}From the theorem it immediately follows
\proclaim{Corollary 1} Let $C$ be as above. For a general line bundle $L$ on
$C$ with the above
lower bounds on the degree, $\gamma_{K,L}$ is onto.
\endproclaim
Furthermore with a little argument one also obtains
\proclaim{Corollary 2} Let $C$ be a general curve
of genus $g>8$ and $L$ be an arbitrary       line bundle on
it, of degree $\ge 5s+12$ if $g=2s$
or $\ge 5s+8$ when $g=2s+1$.
Then $\gamma_{K,L}$ is onto.
\endproclaim

The same attack can be applied to Wahl maps of the kind $\gamma_{L,N}$,
with $L$ and $N$ any two line bundles on $C$. In this direction I prove
the following
\proclaim{Theorem B} Let $C$ be a general curve
of genus $g>8$ and $L$, $N$ be two line
bundles on it. Then
\roster
\item assume that $deg(L)\ge 3s+5$ if $g=2s+1$ (resp., $\ge 3s+4$ if
$g=2s$) , and that $deg(N)\ge deg(L)+g-1$. Then if $L$ and $N$ are chosen
general $\gamma_{L,N}$ is onto.
\item if $L$ and $N$ are arbitrary and $deg(L)\ge 5s+7$ if $g=2s+1$
(resp., $\ge 5s+5$ if $g=2s$) and $deg(N)\ge deg(L)+g-1$, then $\gamma_{L,N}$
is onto.
\endroster
\endproclaim

Corollary 2 and Theorem B (2) should be compared with the similar
results obtained in \cite{BEL}.

The paper is organized as follows. In the first part, Voisin's point
of view is applied to the situation at hand. Specifically,
in $\S$1
we explain the relation between
gaussian maps of the type $\gamma_{K,L}$ and section multiplications
$H^0(C,K-A)\otimes H^0(C,L-A)\longrightarrow H^0(C,K+L-2A)$, for $A$ a pencil
on $C$.
 In particular, it is shown that the proof of the first
statement of Theorem A follows from from the surjectivity of
$$\phi: \bigoplus_{A\in W^1_{s+2}(C)}H^0(K+L-2A)\otimes
\wedge^2H^0(A)\longrightarrow
H^0(2K+L)$$
given by the composition of $id\otimes \gamma_{A,A}$ with section
multiplication. The surjectivity of $\phi$ is dealt with in $\S$2.


In $\S$3 a degeneration argument
is used to show that the above multiplications
are surjective on the general curve, thereby obtaining a surjectivity
statement for $\gamma_{K,L}$ under suitable conditions on $L$.
The proof is given by an induction on the genus.

In $\S$4 these results are extended to the case of the Gaussian maps
$\gamma_{L,N}$, and in $\S$5 an application to higher Wahl maps
is given.

\subheading{Acknowledgments} I wish to express my thanks to my advisor,
Robert Lazarsfeld, for his patient guidance and many very helpful remarks,
and for pointing out several mistakes in the
first draft of this paper, and to Jonathan Wahl for
providing me with useful information about the deformation theoretic
aspect of the problem.

\heading $\S$1. A Basic Commutative Diagram
\endheading

As before, consider a smooth projective curve $C$ and a line bundle
$L$ on it. If $A$ is any other line bundle on $C$ and $V\subset
H^0(A)$ is a pencil of sections of $A$, then we have a commutative diagram

$$
\CD
H^0(L-A)\otimes H^0(K-A)\otimes \wedge^2V @>\alpha>> H^0(L+K-2A)\otimes
\wedge^2(V) \\
\beta @VVV                                 \delta @VVV \\
R(L,K)   @>\gamma_{L,K}>>                   H^0(2K+L)
\endCD \tag {1.1}
$$
where $\alpha$ is section multiplication, and $\beta((s\otimes t)\otimes
(u_1\wedge u_2))=(s\cdot u_1)\otimes (t\cdot u_2)-(s\cdot u_2)\otimes
(t\cdot u_1)$, while $\delta$ is the composition of $id\otimes \gamma_A$
with section multiplication (observe that $\wedge^2(A)\subset \wedge^2(H^0(A))
\subset R(A,A)$).

Recall the following facts about the varieties $W^r_d$:
\roster
\item $W^r_d$ is a connected subvariety of $Pic^d(C)$
\item for any curve $C$, $dim(W^r_d)$ is at least equal to the expected
dimension expressed by the Brill-Noether number $\varrho(r,d,g)=g-(r+1)\cdot
(g-d+r)$
\item if $C$ is Petri general, then $dim(W^r_d)=\varrho (r,d,g)$ and
$W^r_d$ is smooth off $W^{r+1}_d$.
\endroster
In particular it follows that, if $C$ is general and $\varrho (r+1,d,g)<0$,
then $W^r_d$ is a smooth irreducible variety of the expected dimension.

Let us assume in what follows that $C$ is a Petri general curve. Then
$dim(W^1_{s+2})=2$ when $g=2s$ and $=1$ when $g=2s+1$; by the above in both
cases $W^1_{s+2}$ is smooth and irreducible, and nondegenerate in $Pic^d(C)$.
For any nonempty open subset $U\subset W^1_{s+2}$ let us define
$$R_U=:\bigoplus_{A\in U} H^0(L-A)\otimes H^0(K-A)\otimes W_A$$
$$S_U=:\bigoplus_{A\in U} H^0(L+K-2A)\otimes W_A $$
where $W_A=\wedge ^2 H^0(C,A)$.
We then have a commutative diagram
$$
\CD
R_U @>>>            S_U \\
@VVV               \phi @VVV  \\
R(L,K) @>\gamma_{L,K}>>         H^0(2K+L)
\endCD \tag{1.2}
$$
and we clearly have:
\subheading {Observation}
Assume that $\phi$ above is surjective. Then if $\gamma_{L,N}$ is not
onto, not all multiplications $H^0(L-A)\otimes H^0(K-A)\longrightarrow
H^0(L+K-2A)$, for $A\in W^1_{s+2}$, can be surjective.

We can say more:
\proclaim{Proposition 1.1} Let $C$ be a Petri general curve
of genus $g>8$ and
$L$ be a line bundle on it. Suppose that $\phi$ is surjective.
\roster
\item Assume that $deg(L)> 3s$ if $g=2s$,
or that $deg(L)>3s+3$ if $g=2s+1$, and that $L\ne K$ is general.
Define
$U=:\{A\in W^1_{s+2}| h^0(K+A-L)=0\}$. Then if $\gamma_{K,L}$ is not onto
the multiplications $H^0(L-A)\otimes H^0(K-A)\longrightarrow H^0(L+K-2A)$
are never surjective, for any $A\in U$
\item If $\gamma_K$ is not surjective, then $S^2H^0(K-A)\longrightarrow
H^0(2K-2A)$ is not surjective, for any $A\in W^1_{s+2}$
\endroster
\endproclaim
\subheading {Remark}
In (1), the condition on $deg(L)$ implies that $deg(K+A-L)<g$ and so for
a general choice of such an $L$ $U$ will be a nonempty open subset of
$W^1_{s+2}$; the complement of $U$ consists of those points at which the
dimension of $H^0(C,L-A)$ jumps up. In (2), $H^0(K-A)$ has constant rank on
$W^1_{s+2}$. Also, (2) is the content of \cite{V, Lemma 10}.

\demo{Proof} For (1), for a general choice $(A_1,\cdots,A_n)\in U\times
\cdots \times U$ where $n$ is sufficiently large
$$ \bigoplus_{i=1}^nH^0(K+L-2A_i)\longrightarrow H^0(2K+L)$$
is onto.
On the other hand, since $H^0(2K+L)$ and $H^0(K+L-2A)$ have constant rank
on $U$, if one of the above multiplications is onto for some $A\in U$ the
same is true for the general point of $U$.
So if the statement was false a general choice of $(A_1,\cdots,A_n)$ would
yield a composition of surjections
$$
\bigoplus_{i=1}^n H^0(L-A_i)\otimes H^0(K-A_i)\rightarrow
\bigoplus_{i=1}^n H^0(L+K-2A_i)\rightarrow H^0(2K+L)$$
and then $\gamma_{K,L}$ would be onto. For (2), use the same argument
with $U=W^1_{s+2}$.
\enddemo

\input amstex
\documentstyle{amsppt}
\magnification=1200
\hsize=16truecm
\tolerance=2000
\vsize=21truecm
\baselineskip 14pt
\document
\heading $\S 2$. Surjectivity Of $\phi$
\endheading
As we have seen, we are led to the question of the surjectivity of
$$ \phi:\bigoplus_{A\in U} H^0(K+L-2A)\longrightarrow
H^0(2K+L)\tag2--1.$$
We'll prove:
\proclaim{Theorem 2.1} Suppose that $C$ is a Petri general curve
of genus $g>8$ and that $L$ is a general line bundle on $C$, with $deg(L)>
2deg(A)$, and let $U\subset W^1_{2s+2}$ be open and nonempty.
Then $\phi$ is onto.
\endproclaim

\demo{Proof}The following argument is an adaptation of Voisin's. We'll study
thekernel of the dual map
$$\phi^{*}:H^1(T_C-L)\longrightarrow \bigoplus_{A\in U}H^1(2A-L)\tag2--2$$
given by $\phi^{*}(u)=\oplus_{A\in W^1_{s+2}} u\cdot R_A$, where
$R_A\in H^0(2A+K)$ denotes the ramification divisor of the morphism
$\phi_A:C\longrightarrow \bold P^1$ associated to the pencil $A$.
Choose a double cover $\pi:\tilde C\longrightarrow C$ ramified along a general
$B\in |2L|$. We then have the basic isomorphisms
$$ \pi_{*}O_{\tilde C}\simeq O_C\oplus L^{-1}\tag 2--3$$
$$K_{\tilde C}\simeq \pi^{*}(K_C\otimes L)\tag 2--4$$
and so
$$H^1(\tilde C,T_{\tilde C})\simeq H^1(C,T_C-L)\oplus H^1(C,T_C-2L)\tag2--5$$
Hence we can interpret $u\in H^1(C,T_C-L)$ as a first order deformation of
$\tilde C$. It is easily checked that the natural cup product map
$$H^1(C,T_C-L)\longrightarrow Hom(H^0(C,K_C+L),H^1(C,O_C))\tag2--6$$
is a component of
the period map of $\tilde C$; it is furthermore still injective, because its
dual is given by section multiplication $H^0(K_C+L)\otimes H^0(K_C)
\longrightarrow H^0(2K_C+L)$ and this is surjective because $deg(2K+L)>4g+2$
and by a theorem of Mark Green.
So to prove the theorem it is sufficient to show that any $u\in Ker(\phi^{*})$
maps to zero under (1.6).

Observe that, since $L$ is non-torsion, the pull-back map $\pi^{*}:
Pic^{s+2}(C)\longrightarrow Pic^{2s+4}(\tilde C)$ is injective. Hence
$\pi^{*}(W^1_{s+2}(C)\subset Pic^{2s+4}(\tilde C)$ is a smooth subvariety,
isomorphic to $W^1_{s+2}(C)$.
\proclaim{Lemma 2.2}$\pi^{*}(W^1_{s+2})$ is an irreducible component of
$W^1_{2s+4}(\tilde C)$.
\endproclaim
\demo{Proof of Lemma 2.2} Clearly $\pi^{*}(W^1_{s+2}(C))\subset W^1_{2s+4}
(\tilde C)$. By Brill-Noether theory \cite{ACGH} the statement will follow
if we show that, at the general point of $\pi^{*}(W^1_{s+2}(C)$ the Petri
homomorphism $$\mu_{\pi^{*}(A)}:H^0(\tilde C,\pi^{*}(A))\otimes H^0(\tilde
C,K_{\tilde C}\otimes \pi^{*}(-A))\longrightarrow H^0(\tilde C,K_{\tilde C})
\tag2--7$$
has corank equal to the dimension of $W^1_{s+2}(C)$. Assume first that
$A$ is globally generated.
{}From (2.3) and (2.4) one has that (2.7) splits as the direct sum of
$$\mu_A:H^0(C,A)\otimes H^0(C,K-A)\longrightarrow H^0(K_C)\tag2--8$$
$$\nu_A:H^0(C,A)\otimes H^0(C,K+L-A)\longrightarrow H^0(C,K_C+L)\tag2--9$$
and (2.8) is the Petri homomorphism of $A$, which by by the assumption
on $C$ has corank equal to $dim(W^1_{s+2}(C)$, while the base point free
pencil trick applied to $A$, together with the fact that $H^1(K+L-2A)=0$,
show that (2.9) is a surjection.
So the proof will be complete if we show that the general
point $A\in W^1_{s+2}(C)$ is a globally spanned line bundle on $C$.
Consider first the case $g=2s$. Then $dim(W^1_{s+1}(C))=0$ and each $B\in
W^1_{s+1}(C)$ is spanned, because $\varrho(1,d,2s)<0$ $\forall d<s+1$.
Also, if $P\in C$ and $B\in W^1_{s+1}(C)$ then it is easily checked that
$B+P\in W^1_{s+2}(C)$, and that the base point locus of $B+P$ is $\{P\}$.
Hence we get a finite family of (disjoint) copies of $C$ in $W^1_{s+2}(C)$,
one for each element of $W^1_{s+1}$. If $A\in W^1_{s+2}(C)$ is not globally
generated and $P$ is a base
point of $A$, then $B=A-P\in W^1_{s+1}(C)$,
i.e. $A$ lies on one of these curves. Hence it is clear that the lemma holds
in this case. As to the case $g=2s+1$, it is easy to see that all $A\in W^1_
{s+2}(C)$ are globally generated.
\enddemo

Let's now return to the proof of theorem 2.1. Suppose $u\in Ker(\phi^{*})$,
so that for all $A\in U$ we have $A\cdot R_A\in H^1(2A-L)$.
Since $H^0(2A-L)=0$, $H^1(2A-L)$ has constant rank on $W^1_{s+2}(C)$, and so
we actually have $u\cdot R_A=0$ $\forall A\in W^1_{s+2}(C)$.
This has the following deformation theoretic interpretation. First of all,
observe that the first order deformation of $\tilde C$ induced by $u$
carries along a first order deformationof $Pic^{2s+4}(\tilde C)$. Next we
have:
\proclaim{Claim} Suppose that $u\in Ker(\phi^{*})$, and that
$A\in W^1_{s+2}(C)$.
  Then $\pi^{*}(A)$ deforms together with its sections along
the first order deformation induced by $u$.
\endproclaim
\demo{Proof of the Claim} By Brill-Noether theory, a line bundle $N$ on $
\tilde C$
deforms together with its sections along the first order deformation induced
by $\xi \in H^1(\tilde C, T_{\tilde C})$ if and only if $\xi$ annihilates the
image of the Wahl map $\gamma_{N,K_{\tilde C}-N}$ which maps the kernel of
the Petri homomorphism of $N$ to $H^0(\tilde C,2K_{\tilde C})$.
Suppose first that $A$ is base point free.
Applying the base point
free pencil trick to $A$ we get that
$$Ker\mu(\pi^{*}(A))\simeq H^0(C,K+L-2A).$$

Now $\gamma_{L,K_{\tilde C}-L}$ splits as the direct sum of maps
$\alpha_1: H^0(C,K+L-2A)\longrightarrow H^0(C,2K+2L)$ and
$\alpha_2: H^0(C,K+L-2A)\longrightarrow H^0(C,2K+L)$, where
$\alpha_2$ is given by multiplication with $R_A$.
Therefore, since $u\cdot R_A=0$ $u$ annihilates the image of $\alpha_2$.
Because $H^1(K+L)=0$, it also kills the imge of $\alpha_1$, and the Claim
follows in this case.
If $g=2s$ and $A\in W^1_{s+2}$ is not spanned then $A$ has exactly one
base point $P$, and $A-P\in W^1_{s+1}$. Now we apply the base point free
pencil trick to $A-P$, and this yields an exact sequence
$0@>>>H^0(C,K_C+L+P-2A)@>>>H^0(C,A)\otimes H^0(K+L-A) @>>> H^0(C,K+L-P)$
and since the latter space injects into $H^0(C,K+L)$ we obtain
$$Ker(\mu (\pi^{*}A))\simeq H^0(C,K_C+L+P-2A)\subset H^0(C,K_C+L+2P-2A).$$
Now $H^0(C,K_C+2P-2A)\longrightarrow H^0(C,2K_C+L)$ is given by cupping
with $R_A$, and so the statement follows in this case also.
\enddemo

This proves the following:
\proclaim{Corollary 2.3} The first order deformation of $Pic^{2s+4}(C)$
associated to
$u\in Ker(\phi^{*})$ contains a first order deformation of $W^1_{s+2}(C)$.
\endproclaim

Observe that given an inclusion of algebraic manifolds $Y\subset X$
and a first order deformation of $X$ containing a first order deformation
of $Y$ there is a commutative diagram
$$
\CD
H^0(X,\Omega^1_X) @>u_X>> H^1(X,\Cal O_X) \\
@VVV                        @VVV  \\
H^0(Y,\Omega^1_Y) @>u_Y>>   H^1(Y,\Cal O_Y)
\endCD \tag2-10
$$
where $u_X\in H^1(X,T_X)$ and $u_y\in H^1(Y,T_Y)$ are the
extension classes of the two first order deformation.

On the other hand, we have the isomorphisms
$$H^0(X,\Omega^1_{Pic^{2s+4}(\tilde C)})\simeq H^0(C,K+L)\oplus
H^0(C,K)\tag2-11$$
$$H^1(Pic^{2s+4}(\tilde C),\Cal O_{Pic^{2s+4}(\tilde C)})\simeq
H^1(C,\Cal O_C)\oplus H^1(C,-L)\tag 2.12$$
and so we get the following commutative diagram

$$
\CD
H^0(C,K+L) @>u>>         H^1(C,\Cal O_C) \\
a @VVV                      d @VVV \\
H^0(Pic^{2s+4}(\tilde C),\Omega_{Pic^{2s+4}(\tilde C)})  @>>> H^1(Pic^{2s+4}
(\tilde C),\Cal O_{Pic^{2s+4}(\tilde C)})\\
b @VVV                   e @VVV \\
H^0(W^1_{s+2}(C),\Omega_{W^1_{s+2}(C)})  @>>>   H^1(W^1_{s+2}(C),
\Cal O_{W^1_{s+2}(C)}) \\
\endCD \tag 2-13
$$

The proof of Theorem 2.1 is then completed by the following
\proclaim{Lemma 2.4} In diagram (2.13) we have $ba=0$ and $ed$ is
injective.
\endproclaim
\demo{Proof of Lemma 2.4} $b$ is the composition
$$H^0(K_{\tilde C})\longrightarrow H^0(K_C)\hookrightarrow H^0(\Omega^1_W)
\tag 2-14$$
Hence the first assertion follows, because the first map above is just
projection along
$H^0(K+L)$.
Next observe that $ed$ is the composition
$$ H^1(\Cal O_C)\hookrightarrow H^1(\Cal O_{W^1_{s+2}(C)}) \simeq
H^0(\Omega^1_{W^1_{s+2}(C)})\tag 2-15$$
where the last map is injective by \cite{F-L}.
\enddemo
\enddemo
\subheading{Remark}
Since $H^0(C,K_C-A)$ has constant rank on $W^1_{s+2}$, we may
apply this argument to the case $L=K$.
With respect to the proof in \cite{V},
dealing with $W^1_{s+2}$ rather than $W^1_{s+1}$ in the even genus case
avoids the hypothesis $L=K$ and simplifies the argument. However, this is done
at the numerical cost of dealing with pencils that are only next to minimal
rather than minimal in the case of even genus. In other words, when the above
theorem is applied to the particular case $L=K$ and $g=2s$ we only get that if
$\gamma _K$ is not onto then
$K_C-A$ is not projectively normal, for $A\in W^1_{s+2}$, rather than for
$W^1_{s+1}(C)$. In spite of what I was errouneously claiming in a first draft
of this paper this does not imply the stronger numerical statement that
$K_C-A$ is not projectively normal, for $A\in W^1_{s+1}$.

\input amstex
\documentstyle{amsppt}
\magnification=1200
\hsize=16truecm
\tolerance=2000
\vsize=21truecm
\baselineskip 14pt
\document
\heading $\S3$. Surjectivity Of $\gamma_{K,L}$
\endheading

Referring to diagram 1.2, we have shown that under appropriate
conditions the map $\phi$ is onto. It follows as remarked in
$\S1$ that if $\gamma_{K,L}$ is not onto then no multiplication
$H^0(C,K-A)\otimes H^0(C,L-A)\longrightarrow H^0(C,K+L-2A)$ is
onto, for any $A\in U$. This is the statement of Proposition 1.1.
We'll prove:
\proclaim{Theorem 3.1} Let $g\ge 4$ be either $2s$ or $2s+1$, and let
$(C,L,A)$ be a general choice of a smooth curve of genus $g$, and line
bundles $L$ and $A$ on $C$, with $A\in W^1_{s+2}(C)$ and $deg(L)\ge 3s+10$
if $g=2s$ and $deg(L)\ge 3s+6$ for $g=2s+1$. Then $H^0(C,K+A-L)=0$
and section multiplication $H^0(C,K-A)\otimes H^0(C,L-A)\longrightarrow
H^0(C,K+L-2A)$ is surjective.
\endproclaim

Before proving the Theorem, let's remark that it implies the following
\proclaim{Corollary 3.1} Let $C$ be a general curve of genus $g
>8$, with
$g=2s$ or $g=2s+1$. Then if $L$ is a general line bundle on $C$, with
$deg(L)\ge 3s+10$ when $g=2s$ and $deg(L)\ge 3s+6$ when $g=2s+1$ the Wahl
map $\gamma_{K,L}$ is onto.
\endproclaim

{}From this we may deduce a result already contained in \cite{BEL}:
\proclaim{Corollary 3.2} Let $C$ be a general curve
of genus $g>8$ as before, and let $N$ be a line
bundle on $C$, satisfying $deg(N)\ge 5s+12$ when $g=2s$ and $deg(N)\ge
5s+8$ when $g=2s+1$. Then $\gamma_{K,N}$ is onto.
\endproclaim
\demo{Proof of Corollary 3.2} If $B$ is a general line bundle on
$C$ of degree $g+2$, we may assume that $B$ is spanned and that $\gamma_
{K,N-B}$ is onto, by virtue of the previous corollary. Consider a pencil
of sections $V\subset H^0(C,B)$ which generates $V$. We have a
commutative diagram
$$
\CD
R(K,N-B)\otimes V @>\gamma_{K,N-B}>> H^0(C,2K+N-B)\otimes V \\
@VVV                       \beta @VVV \\
R(K,N)    @>\gamma_{K,N}>> H^0(C,2K+N)
\endCD \tag3-1
$$
By assumption, $\gamma_{K,N-B}$ is onto, and the base point free pencil
trick shows that so is $\beta$. Hence $\gamma_{K,N}$ is also onto.
\enddemo

\demo{Proof Of Theorem 3.1} Let us first consider the case $g=2s+1$.
The statement of the Theorem in this case will follow from the following

\proclaim{Proposition 3.3} Let $C$ be a general curve of genus $g>3$, and
let $L$ be a general line bundle on $C$, with $deg(L)\ge 3s+4$ if $g=2s$
and $deg(L)\ge 3s+6$ if $g=2s+1$. Then there exists $A\in W^1_{d_{min}}(C)$
satisfying $H^0(C,K+A-L)=0$ and such that $H^0(C,K-A)\otimes H^0(C,L-A)
\longrightarrow H^0(C,K+L-2A)$ is onto.
\endproclaim

Recall that $d_{min}=s+1$ when $g=2s+1$ and $d_{min}=s+2$ when $g=2s+1$.
The odd genus case of the Proposition is the same as the odd genus case of
the Theorem. The even genus case and the odd genus case of the Proposition
can be proved simoutaneously with an induction argument.

\demo{Proof Of Proposition 3.3} We proceed by induction on $g$.
To begin with, if $C$ is Petri general of genus $g\ge 4$ then by
Riemann-Roch and Brill-Noether number calculations one easily checks that,
$\forall A\in W^1_{d_{min}}(C)$,
$K-A$ is spanned for $g\ge 4$
and birationally very ample for $g>3$
(and very ample for $g\ge 10$).
So let $C$ be Petri general of genus 4, so that $deg(A)=3$ and $K-A$ is
spanned.
In fact, $K-A\in W^1_3(C)$ and by the base point free pencil trick the
surjectivity of section multiplication follows if we have
$H^1(C,L-K)\simeq H^0(C,2K-L)^{*}=0$. But under the given assumptions
$deg(2K-L)\le 2$ and a general choice of $L$ does the job. On the other
hand $deg(K+A-L)<1$ and so $H^0(C,K+A-L)=0$ can also be arranged.
One deals similarly with the case $g=5$.

So now assume given a general curve $C$, which we'll also assume to be
Petri general, of genus $g=2s\ge 4$ and line bundles $A$ and $L$ on $C$ with
$A\in W^1_{s+1}(C)$ and $deg(L)\ge 3s+4$, satisfying the conclusions of
Proposition 3.3.
We may assume without loss that $L-A$ is very ample, and we know that $K-A$
is spanned and birationally very ample. Therefore we have two nondegenerate
morphisms $\varphi_{K-A}:C\longrightarrow \bold P^{s-1}$ and
$\varphi_{L-A}:C\hookrightarrow \bold P^l$, and hence a product embedding
$\varphi:C\hookrightarrow \bold P^{s-1}
\times \bold P^l$.
Choose points $P,Q\in C$ generally and let $l_1\subset \bold P^{s-1}$ and
$l_2\subset \bold P^l$ meeting $C$ nontangentially at $P$ and $Q$ and at
no other point. Identifying $l_1$ and $l_2$ the product embedding gives
a smooth
rational
curve $\Delta\subset \bold P^{s-1}\times \bold P^l$ meeting $\varphi(C)$
nontangentially at $P$ and $Q$.
Define $C^{\prime}=:C\cup \Delta$; then $C^{\prime}$ is a nodal curve,
of genus $g^{\prime}=2s+1$. The proof of the following lemma will be
given later:
\proclaim{Lemma 3.4} $C^{\prime}$ can be smoothed in $\bold P^{s-1}
\times \bold P^l$.
\endproclaim

Next let
$$A^{\prime}=:K_{C^{\prime}}\otimes \Cal O_{\bold P^{s-1}}(-1) \tag3-2$$
Then $deg(A^{\prime})=deg(A^{\prime}|_C)+deg(A^{\prime}|_{\Delta})$
and so
$$deg(A^{\prime})=2g^{\prime}-2-deg(K_C)+deg(A)-1=s+2\tag3-3$$
and by Riemann-Roch
$$h^0(C^{\prime},A^{\prime})=h^0(C^{\prime},K_{C^{\prime}}-A^{\prime})
+s+2+1-g^{\prime}=$$
$$=h^0(C,K-A)+h^0(\bold P^1,\Cal O_{\bold P^1}(1))+2-s=2\tag3-4$$
Next define $L^{\prime}$ on $C^{\prime}$ by
$$L^{\prime}-A^{\prime}=\Cal O_{\bold P^l}(1)|_{C^{\prime}}\tag3-5$$
so that
$$deg(L^{\prime})=deg(A^{\prime})+deg(L-A)+1=deg(L)+2\tag3-6$$
Let's first check that $H^0(C^{\prime},K_{C^{\prime}}+A^{\prime}-L^{\prime})
=0$. This is equivalent to $h^0(C^{\prime},L^{\prime}-A^{\prime})=
deg(L^{\prime}-A^{\prime})+1-g^{\prime}$. But the right hand side is
$h^0(C,L-A)+h^0(\bold P^1,\Cal O_{\bold P^1}(1))-2=h^0(C,L-A)$, while the left
hand side is $deg(L-A)+1-g=h^0(L-A)$.
So this step of the induction is reduced to the following:
\proclaim{Lemma 3.5} $H^0(C^{\prime},L^{\prime}-A^{\prime})\otimes
H^0(C^{\prime},K^{\prime}-A^{\prime})\longrightarrow
H^0(C^{\prime},K^{\prime}+L^{\prime}-2A^{\prime})$ is onto.
\endproclaim
Let's postpone the proof of the above lemma and proceed to the second part
of the induction. So assume given a triple $(C,L,A)$ with $C$ a Petri
general curve of genus $g=2s+1$ and $L$, $A$ line bundles on $C$ satisfying
$A\in W^1_{s+2}(C)$, $deg(L)\ge 3s+6$, $H^0(C,K+A-L)=0$ and such that the
surjectivity statement of the proposition holds.
As before we may assume that $L-A$ is very ample, and we consider the
product embedding $\varphi_{K-A}\times \varphi_{L-A}:C\longrightarrow
\bold P^{s-1}\times \bold P^l$; however, we now consider a copy $l_1\subset
\bold P^{s-1}$ of $\bold P^1$ embedded in degree two, and meeting
$\varphi_{K-A}(C)$ nontangentially at $\varphi_{K-A}(P)$ and
$\varphi_{K-A}(Q)$, $P,Q\in C$, and at no other point. Also, let
$l_2\subset \bold P^l$ be a line meeting $C$ nontangentially at
$\varphi_{L-A}(P)$ and $\varphi_{L-A}(Q)$ and at no other point.
Again, we identify these two copies of $\bold P^1$ and call $\Delta$
the image under the product embedding in $\bold P^{s-1}\times \bold
P^l$. Let $\varphi=:\varphi_{K-A}\times \varphi_{L-A}:C\hookrightarrow
\bold P^{s-1}\times \bold P^l$ and define $C^{\prime}=:C\cup \Delta$.
Then $C^{\prime}$ is a nodal curve of genus $g^{\prime}=g+1$, and
on it we consider the line bundle $A^{\prime}=:K_{C^{\prime}}
\otimes \Cal O_{\bold P^{s-1}}(-1)|_{C^{\prime}}$. By the same computation
as above one checks that $deg(A^{\prime})=s+2$ and $h^0(C^{\prime},
A^{\prime})=2$. Defining $L^{\prime}$ by $L^{\prime}-A^{\prime}=
:\Cal O_{\bold P^l}(-1)|_{C^{\prime}}$ we have $deg(L^{\prime})=
deg(L)+1$. Exactly as before one checks that $dim|L^{\prime}-
A^{\prime}|$ has the expected value, and so the induction step will
be completed by providing a proof of the corresponding variants
of lemma 3.4 and 3.5. Furthermore, observe that in passing from one
odd genus to the next $deg(L)$ increases by three, and this is exactly
what it takes to keep the induction going.
We now prove lemmas 3.4 and 3.5.

\demo{Proof Of Lemma 3.4} Let $Y=\bold P^{s-1}\times \bold P^l$.
The smoothability of $C^{\prime}$ in $Y$ will follow if we can show that
$H^1(C^{\prime},N_{C^{\prime}/Y})=0$ and that $H^0(C^{\prime},
N_{C^{\prime}/Y})\longrightarrow T^1_{C^{\prime}}$ is onto. In fact,
since $C^{\prime}\subset Y$ is a local complete intersection the first
condition implies that the Hilbert scheme of $Y$ is smooth at
$C^{\prime}$, while the second says that there are embedded first order
deformations of $C^{\prime}$ which smooth the nodes (cfr \cite{H}).

Recall that $\varphi_{K-A}$ is
birationally very ample, while $\varphi_{L-A}$ is very ample.
On the other hand, just by Petri generality we have
$H^1(C,T_{\bold P}|_C)=0$ in both cases. From this it is clear
that
$$H^1(C,T_Y|_C)=0\tag 3-8$$

Next recall the exact sequence
$$0@>>>\Cal O_{C^{\prime}}@>>>
\Cal O_C\oplus \Cal O_{\Delta}@>>>\Cal O_{C\cap \Delta}@>>>0\tag3-9$$
from which we obtain the sequence
$$0@>>>T_Y|_{C^\prime}@>>>
T_Y|_C\oplus T_Y|_{\Delta}@>>>T_Y(P)\oplus T_Y(Q)@>>>0\tag3-10$$
where $P$ and $Q$ are the intersection points of $C$ and $\Delta$.
Since $T_Y$ is spanned, the latter sequence is exact on global
sections, and so we get that $H^1(C^{\prime},T_Y|_{C^{\prime}})=0$.
Now the exact sequence
$$0@>>>T_{C^{\prime}}@>>>T_Y|_{C^{\prime}}@>>>N_{C^{\prime}/Y}
@>>>T^1_{C^{\prime}}@>>>0\tag3-11$$
can be chopped off in two short exact sequences, and from this
we see that both of the above conditions are satisfied.
\enddemo

\demo{Proof Of Lemma 3.5}
By assumption, on $C$ we have the exact sequence
$$0@>>>R(K-A,L-A)@>>>H^0(K-A)\otimes H^0(L-A)@>>>H^0(K+L-2A)@>>>0
\tag 3-12$$
where $R(K-A,L-A)\subset H^0(Y,\Cal O_Y)$ is the linear series of the
(1,1)-divisors containing $\varphi(C)$. We have a similar sequence on
$C^{\prime}$:
$$0@>>>R(K^{\prime}-A,L^{\prime}-A^{\prime})@>>>
H^0(L^{\prime}-A^{\prime})\otimes H^0(K^{\prime}-A^{\prime})@>>>
H^0(K^{\prime}+L^{\prime}-2A^{\prime})\tag 3-13$$
and one easily checks that $h^0(L-A)=h^0(L^{\prime}-A^{\prime})$,
$h^0(K-A)=h^0(K^{\prime}-A^{\prime})$ and $h^0(K^{\prime}+L^{\prime}-
2A^{\prime})=h^0(K+L-2A)+1$.
Hence to prove the statement it is sufficient to show that
$dimR(K^{\prime}-A^{\prime},L^{\prime}-A^{\prime})<dimR(K-A,L-A)$,
i.e. that there exist (1,1) divisors in $Y$ containing $C$ but not $\Delta$.

Let us first consider the case $g=2s$, so that $\Delta$ has bidegree
(1,1).
 Observe that then $\Delta \subset Y$ depends on the identification of
of the lines $l_1$ and $l_2$. When we change this identification by an
automorphism of $\bold P^1$, the image of $\Delta$ sweeps the surface
$$<\varphi_{K-A}(P),\varphi_{K-A}(Q)>\times <\varphi_{L-A}(P),
\varphi_{L-A}(Q)>$$ where $<,>$ denotes the line joining the given points
in the appropiate spaces.
The statement in this case follows from the following
\proclaim{Claim} Let $C$ be a smooth curve, and let $\varphi_1:
C\longrightarrow \bold P^m$ and $\varphi_2:C\longrightarrow \bold
P^n$ be nondegenerate morphisms. Define
$$\tilde Sec(C)=:\bigcup_{P,Q\in C}<\varphi_1(P),\varphi_1(Q)>
\times <\varphi_2(P),\varphi_2(Q)>\tag3-14$$
Then $\tilde Sec(C)$ is not contained in any (1,1)-divisor $D\subset
\bold P^m\times \bold P^n$.
\endproclaim
\demo{Proof Of The Claim} Provisionally let $Y=\bold P^m\times \bold
P^n$ and for $P\in C$ fixed consider the projections $\pi_1:
\bold P^m\setminus \{\varphi_1(P)\}\longrightarrow \bold P^{m-1}$ and
$\pi_2:\bold P^n\setminus \{\varphi_2(P)\}\longrightarrow \bold P^{n-1}$,
where $\bold P^{m-1}$ and $\bold P^{n-1}$ are two fixed hyperplanes.
We then have a product morphism $$\pi_P:(\bold P^m\setminus
\{\varphi_1(P)\})\times
(\bold P^n \setminus \{\varphi_2(Q)\})\longrightarrow \bold P^{m-1}\times
\bold P^{n-1}.$$
Let $D^{\prime}=:D\cap \bold P^{m-1}\times \bold P^{n-1}$,
$\tilde C=:\pi_P(C)$. I claim that
$$D\supset \tilde Sec(C) \Rightarrow D^{\prime}\supset
\tilde Sec(\tilde C)\tag 3-15$$
Note that by induction this reduces the proof of the Lemma to the case where
either $m=1$ or $n=1$, and then it is trivial.

In fact it is easily checked that $\pi_P(\tilde Sec(C))\subset
\tilde Sec(\tilde C))$ is a Zariski dense open subset, so that a general
point in the latter variety can be written as $\pi_P(Q)$, where
$Q=(Q_1,Q_2)\in <\varphi_1(A),\varphi_1(B)>\times <\varphi_2(A),
\varphi_2(B)>$ for suitable $A,B\in C$.
Assume that in this situation
$$D\supset <\varphi_1(P),Q_1>\times <\varphi_2(P),Q_2>\tag 3-16$$
It then follows that
$$D^{\prime}\supset [<\varphi_1(P),Q_1>\times <\varphi_2(P),Q_2>]\cap
(\bold P^{m-1}\times \bold P^{n-1})$$
and this establishes (3-15).

It remains to prove 3-16. By assumption we have $D\supset [<\varphi_1
(P),\varphi_1(A)>\times <\varphi_2(P),\varphi_2(A)>]\cup
[<\varphi_1(P),\varphi_1(B)>\times <\varphi_2(P),\varphi_2(B)>]
\times [<\varphi_1(A),\varphi_1(B)>\times <\varphi_2(A),\varphi_2(B)>]$.
In other words, we have $$D\supset (l_1\times l_2)\cup (r_1\times
r_2)\cup (s_1\times s_2)$$ where $l_1,r_1,s_1\subset \bold P^m$ and
$l_2,r_2,s_2\subset \bold P^n$ are triples of lines contained in the same
plane $\Lambda_1$ and $\Lambda_2$ respectevely.
Intersecting $D$ with $\Lambda_1\times \Lambda_2$ we are thus reduced
to proving that a (1,1)-divisor in $\bold P^2\times \bold P^2$ cannot
contain such a union of lines, and this is a well known fact.
A simpler version of the same argument deals with the case $g=2s+1$.
\enddemo

Let's now come to the even genus case of Theorem 3.1. Recall that we
want to prove the following statement: if $(C,L,A)$ is a general choice
of a smooth curve of genus $g=2s$ and line bundles $L$ and $A$ on $C$,
with $A\in W^1_{s+2}(C)$ and $deg(L)\ge 3s+10$, then $H^0(C,K+A-L)=0$
and section multiplication $H^0(C,K-A)\otimes H^0(C,L-A)\longrightarrow
H^0(C,K+L-2A)$ is onto.
Clearly the first conclusion is true, and
the second follows from a slight modification of the previous
degeneration argument, as follows. Let $(C,L,A)$ be a general choice of
a curve of odd genus $g=2s+1$ and line bundles $L$ and $A$ on it,
satiafying both conclusions of Theorem 3.1. Such a choice exists by the
argument above. Now apply the same construction, but take $\Delta$
to be of bidegree $(1,1)$. $C^{\prime}$ has genus $g^{\prime}=2s+1$,
and if we define $A^{\prime}$ and $L^{\prime}$ in the same
manner it we obtain
$deg(A^{\prime})=s+3$ and $h^0(C^{\prime},A^{\prime})=2$; everything else
stays unchanged.
\enddemo

One may also obtain the statement for genera divisible by 4
from the following covering argument.
Let $(C^{\prime},L^{\prime},A^{\prime})$ be a general choice
of a curve of genus $g^{\prime}=2k+1$ and of line bundles $L^{\prime}$
and $A^{\prime}$ on it satisfying both conclusions of theorem 3.1. Such
a triple exists by the argument above. Pick a general point $P\in C^{\prime}$;
without loss of generality, we may assume that if $N=:L^{\prime}+3P$,
then $H^0(C^{\prime},N-A^{\prime})\otimes H^0(C^{\prime},K^{\prime}-
A^{\prime})\longrightarrow H^0(C^{\prime},N+K^{\prime}-2A^{\prime})$ is also
onto, where we have let $K^{\prime}=:K_{C^{\prime}}$. Also, we may assume
$C^{\prime}$ to be Petri general, so that $d_{min}=k+2$.

Now consider the double cover $C@>\pi>> C^{\prime}$ ramified along $R=6P$.
Letting $g$ be the genus of $C$, we have $g=2(g^{\prime}+1)$ and so
$d_{min}(C)+1=g^{\prime}+3=2(k+2)$. Now let $A=:\pi^{*}(A^{\prime})$,
so that $deg(A)=2deg(A^{\prime})=2(k+2)$ and $h^0(A)=2$, and let
$L=:\pi^{*}(N)$, so that $deg(L)\ge 3(g^{\prime}+1)+10$. Finally observe
that the multiplication map $H^(C,L-A)\otimes H^0(C,K-A)\longrightarrow
H^0(C,L+K-2A)$ splits as the direct sum of various section multiplications
on $C^{\prime}$, two of which are $H^0(C^{\prime},K^{\prime}-A^{\prime})
\otimes H^0(C^{\prime},L^{\prime}-A^{\prime})\longrightarrow
H^0(C^{\prime},K^{\prime}+L^{\prime}-2A^{\prime})$ and
$H^0(C^{\prime},K^{\prime}-A^{\prime})\otimes H^0(C^{\prime},
N^{\prime}-A^{\prime})\longrightarrow  H^0(C^{\prime}, K^{\prime}+N^{\prime}-
2A^{\prime})$
{}From this it is easy to deduce the statement.
\enddemo
\enddemo

\input amstex
\documentstyle{amsppt}
\magnification=1200
\hsize=16truecm
\tolerance=2000
\vsize=21truecm
\baselineskip 14pt
\document
\heading $\S4$. The Wahl Maps $\gamma_{L,N}$
\endheading
As before, consider a smooth projective curve $C$ and let $L,N$ be
two line bundles on it. Suppose $A\in W^1_{s+2}$.
Then $r(L-A),r(N-A)\ge 0$ when $deg(L),deg(N)\ge 3s+2$ in case $g=2s$
and $deg(L),deg(N)\ge 3s+3$ when $g=2s+1$.
Define
$$U=:\{A\in W^1_{s+2}|h^0(C,K+A-L)=h^0(K+A-L)=0\}\tag 4.1$$
and consider the variant of diagram (1-2), where we let
$W_A=:\wedge^2H^0(C,A)$:
$$
\CD
\bigoplus_{A\in U}H^0(L-A)\otimes H^0(N-A)\otimes W_A@>>> \bigoplus_{A\in U}
H^0(N+L-2A)\otimes W_A \\
@VVV                      \phi  @VVV \\
R(L,N)    @>\gamma_{L,N}>>  H^0(K+L+N)
\endCD \tag 4-2
$$
\proclaim{Proposition 4.1} Let $C$ be a Petri general curve and
suppose $deg(L),deg(N)\ge 3s+2$ if $g=2s$ and $deg(L),deg(N)\ge 3s+3$
if $g=2s+1$.
Then $\phi$ from diagram (4-2) is onto.
\endproclaim
\demo{Proof Of Proposition 4.1} Write $N=K-R$, where $deg(R)\le s-4$
if $g=2s$, $deg(R)\le s-3$ if $g=2s+1$.
We look at
$$\phi^{*}:H^1(C,T_C-(L-R))\longrightarrow
 \bigoplus _{A\in U}H^1(C,2A-(L-R))\tag 4-3$$
Consider the double cover $\pi:\tilde C@>>>C$ ramified along the general
element $B\in |2(L-R)|$. The argument used in the proof of Theorem 2.1
goes over verbatim.
\enddemo

\proclaim{Corollary 4.2} In the situation of the Proposition,
if $\gamma_{L,N}$ is not onto then no multiplication
$H^0(C,L-A)\otimes H^0(C,N-A)\longrightarrow H^0(C,L+N-2A)$ is onto,
for any $A\in U$.
\endproclaim

We can now prove
\proclaim{Theorem 4.3} Let $C$ be a general curve and $L,N$ be two
general  line
bundles on it. Assume that
\roster
\item $deg(L)\ge 3s+4 $ if $g=2s$, or that $deg(L)\ge 3s+5$ if
$g=2s+1$
\item $deg(N)>deg(L)+g-2$
\endroster
Then $\gamma_{L,N}$ is onto. Furthermore, assume that $L,N$ are
any two line bundles on $C$ such that (2) holds and $deg(L)\ge 5s+7$
if $g=2s$ or $deg(L)\ge 5s+5$ when $g=2s+1$. Then $\gamma_{L,N}$ is onto.
\endproclaim
\demo{Proof} Let's start with the first statement.
For a general choice of such an $L$ and a general $A\in
U$ $L-A$ is spanned. Let $V\subset H^0(C,L-A)$ be a pencil spanning
$L-A$, so that there is an exact sequence $0@>>>A-L@>>>V\otimes
\Cal O_C@>>>L-A@>>>0$. Twisting by $N-A$ we get the exact sequence
$$0@>>>N-L@>>>V\otimes (N-A)@>>>L+N-2A@>>>0\tag 4-4$$
and this exact on global sections because by condition (2) and the
assumed generality of $L$ and $N$ we may assume $H^1(C,N-L)=0$.
The first statement follows.

To prove the second part of the theorem, observe that for general line
bundles $B$ and $D$ on $C$ of degree $g+1$ we may assume that $B$ and $D$
are spanned and that $\gamma_{N-D,L-B}$ is onto, by the first
statement of this theorem.
We may also assume that $r(B)=r(D)=1$, so that we have exact sequences
$0@>>>-B@>>>H^0(C,B)@>>>B@>>>0$ and $0@>>>-D@>>>H^0(C,D)@>>>D@>>>0$.
Letting $V=:H^0(C,B)$, $W=:H^0(C,D)$ we obtain the commutative diagram
$$
\CD
R(L-B,N-C)\otimes V\otimes W @>\gamma_{L-B,N-D}>> H^0(C,K+L-B+N-C)\otimes
V\otimes W \\
@VVV                  \beta   @VVV     \\
R(L,N) @>\gamma_{L,N}>> H^0(C,K+L+N)
\endCD\tag 4-5
$$
and using the above sequences it is easy to see that $\beta$ is onto.
The theorem follows.
\enddemo

\input amstex
\documentstyle{amsppt}
\magnification=1200
\hsize=16truecm
\tolerance=2000
\vsize=21truecm
\baselineskip 14pt
\document
\heading{$\S$5. Higher Wahl Maps}
\endheading

Consider as before a smooth projective curve $C$ and line bundles
$L$ and $N$ on it. The map $\gamma_{L,N}$ that we have been considering
so far generalizes to a hierarchy of Wahl maps
$\gamma^i_{L,N}$ defined as follows (cfr \cite{W1}). On the product
$C\times C$ consider the line bundle $p_1^{*}(L)\otimes p_2^{*}(N)$,
which we'll abbreviate to $L_1\otimes N_2$, and let $\Delta \subset C\times
C$ denote the diagonal. Then $H^0(C\times C,L_1\otimes N_2)\simeq
H^0(C,L)\otimes H^0(C,N)$ is filtered by the subspaces
$R_l(L,N)=H^0(C\times C,L_1\otimes N_2(-l\Delta))$, $l\ge 0$. For each $l$, we
have an exact sequence on $C\times C$
$$0@>>>L_1\otimes N_2(-(l+1)\Delta)@>>>L_1\otimes N_2(-l\Delta)@>>>
L_1\otimes N_2(-l\Delta)|_{\Delta}@>>>0\tag 5-1$$
and we simply define
$$\gamma^l_{L,N}:R_l(L,N)\longrightarrow
H^0(C,lK_C+L+N)\tag5-2$$
to be the induced map on global sections. Notice that $R_{l+1}(L,N)=
Ker\gamma^l_{L,N}$. $\gamma^0_{L,N}$ is just section multiplication, and
$\gamma^1_{L,N}$ is the usual Wahl map.
The approach used in the previous paragraphs to deal with the first Wahl map
can be generalized to $\gamma^l_{L,N}$, as follows.
As before, let's assume the genus of $C$ is $g=2s$ or $g=2s+1$, and let's
denote $W=W^1_{s+2}$.
Let $U=:\{(A_1,\cdots,A_l)\in W\times \cdots \times W|
h^0(C,K+\sum_{i=1}^lA_i-L)=h^0(C,K+\sum_{i=1}^lA_i-N)=0\}$, and define

$$V(L,N;A_1,\cdots,A_l)=:H^0(L-\sum_{i=1}^lA_i)\otimes
H^0(N-\sum_{i=1}^lA_i)\otimes \wedge^2H^0(A_1)\otimes \cdots
$$
and
$$W(L,N;A_1,\cdots,A_l)=:H^0(L+N-2\sum_{i=1}^lA_i)\otimes
\wedge^2H^0(A_1)\otimes \cdots \otimes \wedge^2 H^0(A_l)$$
and consider the following commutative diagram, which generalizes
(4-2):

$$
\CD
\bigoplus_{(A_1,\cdots,A_l)\in U}V(L,N;A_1,\cdots,A_l) @>>>
\bigoplus_{(A_1,\dots,A_l)\in U}W(L,N;A_1,\cdots,A_l) \\
f @VVV                          \phi @VVV \\
R_l(L,N) @>\gamma^l_{L,N}>> H^0(C,lK+L+N)
\endCD \tag5-3
$$

Assume that $\phi$ is onto. Then, by the usual argument, if
$\gamma^l_{L,N}$ is not onto none of the multiplications
$H^0(L-
\sum_{i=1}^lA_i)\otimes H^0(N-\sum_{i=1}^lA_i)\longrightarrow
H^0(L+N-2\sum_{i=1}^lA_i)$ can be surjective.
Now we can proceed inductively to draw the same kind of conclusions
as in the case of $\gamma_{L,N}$.

\proclaim{Theorem 5.1} Assume that $C$ is Petri general and that
$deg(L)+deg(N)>2g-2+2l(s+2)$. Then
$$\phi:\bigoplus_{(A_1,\cdots,A_l)\in U}H^0(L+N-2\sum_{i=1}^nA_i)
\longrightarrow H^0(lK+L+N)$$
is onto.
\endproclaim
\demo{Proof Of Theorem 5.1}
We have
$$\phi^{*}:H^1((l-1)T-L-N)\longrightarrow \bigoplus H^1(2\sum_{i=1}^lA_i+K
-L-N)$$
given by $u\mapsto \oplus u\cdot R_{A_1}\cdots R_{A_l}$.
If $u\in Ker\phi^{*}$, then $u\cdot R_{A_1}\cdots R_{A_l}=0$,
for every $(A_1,\cdots,A_l)\in W$.

Now $u\cdot R_{A_1}\cdots R_{A_{l-1}}\in H^1(T-(L+N-K-2\sum_{i=1}^{l-1}A_i))$,
and we have $deg(L+N-K-\sum_{i=1}^{l-1}2A_i)>2deg(A)$.
By theorem 2.1, this implies $u\cdot R_{A_1}\cdots R_{A_{l-1}}=0$ for all
$(A_1,\cdots,A_{l-1})\in U$, and now the statement follows by induction.
\enddemo

\proclaim{Corollary 5.2} In the situation of Lemma 5.1, if $\gamma^l_{L,N}$
fails to be surjective then $$H^0(C,L-\sum_{i=1}^lA_i)\otimes
H^0(C,L-\sum_{i=1}^lA_i)\longrightarrow H^0(C,L+N-2\sum_{i=1}^lA_i)$$
is never surjective, for any $(A_1,\cdots,A_l)\in U$.
\endproclaim

\proclaim{Theorem 5.3} Assume that $L$ and $N$ are two general line bundles on
$C$, with $deg(L)\ge g+1+l(s+2)$ and $deg(N)>g-2+deg(L)$. Then
$\gamma^l_{L,N}$ is onto. If $L$ and $N$ are arbitrary and $deg(L)>2g+2
+l(s+2)$ and $deg(N)>g-2+deg(L)$, then $\gamma^l_{L,N}$ is onto.
\endproclaim
\demo{Proof Of Theorem 5.3} Use the same argument as in the proof of theorem
4.3, replacing $L-A$ by $L-\sum_{i=1}^lA_i$ and similarly for $N-A$.
\enddemo

\bigpagebreak
\bigpagebreak
\bigpagebreak

\input amstex
\documentstyle{amsppt}
\magnification=1200
\hsize=16truecm
\tolerance=2000
\vsize=21truecm
\baselineskip 14pt
\document
\subheading{Bibliography}

\cite{ACGH} E. Arbarello, M. Cornalba, P. Griffiths, J. Harris
"Geometry Of Algebraic Curves" Vol.1, Springer Verlag (1984)

\cite{BM} Beauville and Merindol "Section Hyperplanes Des
Surfaces K-3" Duke Math. J. 55 (1987), 873-878

\cite{BEL} A. Bertram, L. Ein, and R. Lazarsfeld "Surjectivity Of
Gaussian Maps For Line Bundles Of Large Degree On Curves" preprint

\cite{CHM} C. Ciliberto, J. Harris, H. P. Miranda "On The Surjectivity
Of Wahl Maps" Duke Math. J. 57 (1988) 829-858

\cite{F-L} W. Fulton, R. Lazarsfeld "On The Connectedness Of Degeneracy
Loci And Special Divisors"  Acta Mathematica 146 (1981), 271-283

\cite{H-H} R. Hartshorne, A. Hirschowitz "Smoothing Algebraic Space Curves"
Lect. Notes In Mathematics 1124, Algebraic Geometry Proceedings,
Sitjes (1983)

\cite{M} S. Mukai "Curves, K3 surfaces, and Fano three-folds of genus
$10$" in Algebric Geometry and Commutative Algebra In Honor Of M. Nagata
(1987), 357-377

\cite{T} Tendian "Deformations Of Cones Over Curves Of High Degree"
Ph.D. dissertation (University Of Northern Carolina, Chapel Hill)
July 1990

\cite{V} C. Voisin "Sur l'Application de Wahl des Courbes Satisfaisantes
La Condition De Brill-Noether-Petri" preprint

\cite{W1} J. Wahl "Deformations Of Quasihomogeneous Surface Singularities"
Math. Ann. 280 (1988) 105-128

\cite{W2} J. Wahl "Introduction To Gaussian Maps On An Algebraic Curve"
Notes prepared in connection with lectures at the Trieste conference
on projective varieties, 1989

\cite{W3} J. Wahl "Gaussian Maps On Algebraic Curves" J. Of Diff. Geometry
32 (1990), 77-98

\cite{W4} J. Wahl "Gaussian Maps And Tensor Products Of Irreducible
Representations" preprint

\bigpagebreak

\subheading{Address}

Mathematics Department

UCLA

Los Angeles, CA 90024

\centerline{September 1992}

\enddocument